\documentclass[a4paper, 12pt]{amsart}
\usepackage{amsmath, amssymb, eucal, amscd, amstext}

\usepackage{enumerate}

\newtheorem{theorem}{Theorem}[section]
\newtheorem{corollary}[theorem]{Corollary}
\newtheorem{lemma}[theorem]{Lemma}
\newtheorem{proposition}[theorem]{Proposition}
\newtheorem{conjecture}[theorem]{Conjecture}

\theoremstyle{definition}

\theoremstyle{remark}
\newtheorem{remark}[theorem]{Remark}

\newtheorem{example}[theorem]{Example}
\numberwithin{equation}{section}

\newcommand{\smnoind}{\smallskip\noindent}

\newcommand{\CL}{\mathcal{L}}

\newcommand{\CK}{\mathcal{K}}

\newcommand{\cv}{\mathbf{c}}
\newcommand{\II}{I\!I}
\newcommand{\her}{\operatorname{her}}

\begin{document}

\baselineskip=15pt

\title[Linear orthogonality preservers: the real rank zero case]{Linear orthogonality
preservers of Hilbert $C^*$-modules over $C^*$-algebras with real rank zero}

\author{Chi-Wai Leung, Chi-Keung Ng  \and Ngai-Ching Wong }

\address[Chi-Wai Leung]{Department of Mathematics, The Chinese
University of Hong Kong, Hong Kong.}
\email{cwleung@math.cuhk.edu.hk}

\address[Chi-Keung Ng]{Chern Institute of Mathematics and LPMC, Nankai University, Tianjin 300071, China.}
\email{ckng@nankai.edu.cn}

\address[Ngai-Ching Wong]{Department of Applied Mathematics, National Sun Yat-sen University,  Kaohsiung, 80424, Taiwan, R.O.C.}
\email{wong@math.nsysu.edu.tw}

\thanks{The authors are supported by Hong Kong RGC Research Grant (2160255),
National Natural Science Foundation of China (10771106), and
Taiwan NSC grant (NSC96-2115-M-110-004-MY3).}

\date{\today}

\keywords{Hilbert $C^*$-modules, orthogonality preservers, modules maps, local maps, real rank zero $C^*$-algebras}

\subjclass[2000]{46L08, 46H40}

\begin{abstract}
Let $A$ be a $C^*$-algebra.
Let $E$ and $F$ be Hilbert $A$-modules with $E$ being full.
Suppose that $\theta : E\to F$ is a linear
map preserving orthogonality, i.e.,
$$
\langle \theta(x),\theta(y)\rangle\ =\ 0\quad\text{whenever}\quad\langle x,y\rangle\ =\ 0.
$$
We show in this  article that if, in addition, $A$ has real
rank zero,
and $\theta$ is an $A$-module map (not assumed to be bounded), then there exists a central positive multiplier
$u\in M(A)$ such that
$$
\langle \theta(x), \theta(y)\rangle\ =\ u
\langle x, y\rangle\qquad (x,y\in E).
$$
In the case when $A$ is a
standard $C^*$-algebra or when $A$ is a $W^*$-algebra containing no
finite type $\II$ direct summand, we also get the same conclusion with the
assumption of $\theta$ being an $A$-module map weakened to being a local map.
\end{abstract}

\maketitle

\section{Introduction and Notations}

It is a common knowledge that, together with the linearity,
the inner product and the norm structures of a  Hilbert space $H$
determine each other.
It might be a bit less well-known that the orthogonality
structure also suffices to determine the inner product up to a scalar.
This fact follows from the following easy observation:
for any $x,y\in H$,
$\|x\| = \|y\| \text{ if and only if } x+\lambda y \text{ is orthogonal to } x-\lambda y \text{ for some scalar }
\lambda  \text{ with } |\lambda|=1$ (see also \cite{Blanco06, chmielinski05}).

It is natural and interesting to ask whether the linear structure and orthogonality structure of a (complex) Hilbert $C^*$-module determines its $C^*$-algebra-valued inner product.
More precisely, let $A$ be a (complex) $C^*$-algebra, and $\theta : E\to F$ be a $\mathbb{C}$-linear map between
Hilbert $A$-modules that preserves orthogonality (i.e. preserves zero $A$-valued inner products).
We want to study to what extent,
$\theta$ will respect the $A$-valued inner products. When the underlying $C^*$-algebra is $\mathbb C$, it reduces to
the case of Hilbert spaces.

We first note that without any further assumption on $\theta$, the above  question might have a negative answer.

\begin{example}\label{eg:module-ortho-map}
Let $H$ be an infinite dimensional (complex) Hilbert space and $A=\mathcal{K}(H)$ be the
$C^*$-algebra of all compact operators on $H$. Suppose that $\bar H$ is a
vector space that is conjugate-linear isomorphic to $H$. When
equipped with the operations: $\langle \overline{\eta_1},
\overline{\eta_2} \rangle := \eta_1\otimes\eta_2$ and
$\overline{\eta_1} T := \overline{T^* \eta_1}$ ($\overline{\eta_1},
\overline{\eta_2}\in \bar H; T\in A$), we see that $\bar H$ is a
Hilbert $A$-module. Suppose that $\theta$ is any unbounded bijective
$\mathbb{C}$-linear map from $\bar H$ onto $\bar H$. Since $\langle x, y \rangle =
0$ if and only if $x =0$ or $y = 0$, we see that both $\theta$ and $\theta^{-1}$ preserves orthogonality.
\end{example}

As we are dealing with Hilbert $A$-modules, a natural additional assumption is that $\theta$ is an $A$-module map, i.e., $\theta(xa)=\theta(x)a$ ($x\in E, a\in A$).
In \cite{Turnsek-JMAA},
Ili\v{s}evi\'{c} and Turn\v{s}ek showed that if
$A$ is a standard $C^*$-algebra, then for every
orthogonality preserving
$A$-module map $\theta: E\to F$, there is a scalar $\lambda\geq 0$ such that
$\langle \theta(x),\theta(y) \rangle = \lambda \langle x,y\rangle$ ($x,y\in E$).
In particular, all such $\theta$ are scalar multiples of isometries.

In \cite{LNW-orth-pres-comm}, under a weaker assumption on $\theta$, namely $\theta$ being ``local'',
we get the same conclusion in the case when $A$ is
a commutative $C^*$-algebra (in fact,  the main difficulties in \cite{LNW-orth-pres-comm} come from the fact that $A$ is not assumed to be an $A$-module map).
Recall that a $\mathbb{C}$-linear map $\theta: E\to F$
is \emph{local} if
$$
\theta(x)a = 0 \quad\text{whenever}\quad xa=0\qquad (x\in E;  a\in A ).
$$
Readers should find the idea of local linear maps familiar. For
example,  local linear
maps on the space of smooth functions defined on a manifold modeled
on $\mathbb{R}^n$ are exactly linear differential operators (see,
e.g., \cite{Pe60, Na85}). See also \cite{KN, Ar04} for the vector-valued case, and
\cite{Al08} for the Banach $C^1[0,1]$-module setting.
Notice that every module map is local, but local linear maps, e.g., linear differential operators,
might not be a module map.
Nevertheless, it has been shown in \cite[Proposition A.1]{LNW-auto-cont} that every \emph{bounded} local map between Hilbert $C^*$-modules is a module map.

The results in \cite{Turnsek-JMAA} and \cite{LNW-orth-pres-comm} lead to the following conjecture.
We remark that the fullness assumption of $E$ in this conjecture is a necessity.
Without this, the conclusion does not hold even in the case when $A$ is commutative (see \cite[3.6]{LNW-orth-pres-comm}).
Here, a Hilbert $A$-module $E$ is said to be
\emph{full} if the linear span of
$\{\langle x, y\rangle: x,y\in E\}$
is dense in $A$.

\begin{conjecture}
\label{conj}
Let $A$ be a $C^*$-algebra.
Let $E$ and $F$ be Hilbert $A$-modules with $E$ being full.
If $\theta: E\to F$ is a ($\mathbb C$-linear) local map preserving orthogonality, i.e. for any $x,y\in E$,
$$
 \langle x,y\rangle\ =\ 0 \quad\text{implies}\quad \langle \theta(x),\theta(y)\rangle\ =\ 0,
$$
then there is a central positive element $u\in M(A)$ such that
\begin{align*}
\langle \theta (x), \theta(y) \rangle\ =\ u \langle x, y \rangle \qquad (x,y\in E).
\end{align*}
\end{conjecture}

In this article, positive answers of this conjecture will be given in the following three cases:
\begin{enumerate}
\item $A$ is a $C^*$-algebra of real rank zero and $\theta$ is an $A$-module map (Theorem \ref{thm:unbdd-case-rro}).
\item $A$ is a standard $C^*$-algebra (Corollary \ref{cor:ortho-local-standard}).
\item $A$ is a $W^*$-algebra with no finite type $\II$ direct summand (Corollary \ref{cor:infiniteWalg}).
\end{enumerate}

As a final remark for the introduction, we note that, unlike the situation in some other literatures (e.g. \cite{FMP}),
$\theta$ is \emph{not} assumed to be bounded, for the conceptual reasons as stated in the beginning of this Introduction (but whose boundedness will be an automatic consequence of our results).

Let us now give some notations that will be used throughout this article.
In the following, $A$ is a $C^*$-algebra, $E$ and $F$ are
Hilbert $A$-modules, and $\Psi, \theta: E \rightarrow F$ are orthogonality preserving $\mathbb{C}$-linear maps, which  are not assumed to
be bounded.

Let $a\in A_+$.
We set $C^*(a)$ to be the $C^*$-subalgebra
generated by $a$, and $\cv(a)$ to be the \emph{central cover} of $a$, i.e., the smallest central element in $A^{**}_+$ dominating $a$ (see, e.g., \cite[2.6.2]{Ped}).
If, in addition, $\alpha, \beta\in \mathbb{R}_+$, we
put $e_a(\alpha, \beta]$ to be the spectral
projections of $a$ in $A^{**}$ corresponding to the set $(\alpha, \beta]\cap \sigma(a)$.

We denote by $Z(A)$ the center and by $M(A)$ the
space of all multipliers of $A$.
On the other hand, ${\rm
Proj}_1(A)$ is the set of all norm-one (i.e., non-zero) projections
in $A$.
For any \emph{open projection} $p\in {\rm Proj}_1(A^{**})$ (i.e.,
there exists an increasing net $\{a_i\}$ in $A_+$ such that
$a_i\uparrow p$ in the weak-*-topology), we denote by $\her(p) :=
pA^{**}p\cap A$ the hereditary $C^*$-subalgebra associated to $p$.
See, e.g., \cite{Bro88} for more information about open projections.

\section{Orthogonality preserving $A$-module maps when $A$ has real rank zero}

We recall that $A$ has \emph{real rank zero} if every self-adjoint element in $A$
can be approximated in norm by invertible self-adjoint elements.
This implies, in particular, that the linear span of ${\rm Proj}_1(A)$ is norm dense in $A$ (see, e.g., \cite{Brown91}).

Let us start with the following easy lemma.
Part (a) of which might be well-known but we give an argument here for completeness.

\begin{lemma}\label{lem:cp}
(a) If $p\in {\rm Proj}_1(A^{**})$ and $b\in Z(pA^{**}p)_+$, then $\|\cv(b)\| = \|b\|$, $\cv(b)p = b$ and $\cv(b)\cv(p) = \cv(b)$.

\smnoind
(b) Suppose that $A$ has real rank zero and $E$ is full.
If $q\in A^{**}\setminus \{0\}$ is an open projection, there are $r\in {\rm Proj}_1(A)$ and $y\in Er$ such that $r = \langle  y,y \rangle  \leq q$.
\end{lemma}
\begin{proof}
(a) Since $b\leq \|b\| 1$, we see that $0\leq b\leq \cv(b)\leq \|b\| 1$ and $\|\cv(b)\| = \|b\|$.
Clearly, $\cv(b)p = p\cv(b)p \geq pbp = b$.
Conversely, as $Z(pA^{**}p) = Z(A^{**})p$ (see e.g. \cite[5.5.6]{KR}), there is $a\in Z(A^{**})_+$ with $b = ap$ (note that $b^{1/2}\in Z(A^{**})p$).
Thus, we have
$b = a^{1/2}pa^{1/2} \leq a^{1/2}\cv(p)a^{1/2} = a\cv(p)$.
As $a\cv(p)$ is central, $\cv(b)\leq a\cv(p)$ and $\cv(b)p = p\cv(b)p \leq ap = b$.
The last equality follows from \cite[2.6.4]{Ped} and the fact that $b\cv(p) = b$.

\smnoind (b) Note that $\her(q)\neq (0)$, and also has real rank zero (see e.g. \cite{Brown91}).
Moreover, $E_0 := E\cdot \her(q)$ is a full (and, hence non-zero) Hilbert $\her(q)$-module.
Pick any $x\in E_0$ such that $a := \langle x,x\rangle $ is a norm one
element in $\her(q)$.
Let $t\in (0,1/3)$ and $b\in \her(q)_+$ such
that $\|a - b\| < t$ and $\sigma(b) = \{\lambda_1,...,\lambda_n\}$
with $\lambda_1 \leq \cdots \leq \lambda_n = \|b\|$ (see e.g.
\cite{Brown91}).
Since $\|b\|>2/3$,  we can choose $s\in
[t,\|b\|]\setminus \sigma(b)$.
If we set $r := e_b(s,2]\in {\rm Proj}_1(A)$, then $\|r -
rar\| \leq \|r - rbr\| + \|b -a\| < 1$.
If $c := r +
\sum_{n=1}^\infty (r-rar)^n\in A_+$, then $(rar)c=c(rar)=r$ and so,
$\langle  xc^{1/2}, xc^{1/2}\rangle = c^{1/2}rarc^{1/2} = r$.
Finally, $xc^{1/2}\in Er$
as $c^{1/2}$ is in the $C^*$-subalgebra $rAr + \mathbb{C}r$.
\end{proof}

\begin{proposition}\label{prop:bddcase-with1}
Let $A$ be a unital $C^*$-algebra of real rank zero.
Suppose that $\theta:E\to F$ is an $A$-module map preserving orthogonality, and there is an element $x_0\in E$ such
that $\langle x_0,x_0\rangle =1$.
Then one can find $u\in Z(A)_+$ satisfying
$$
\langle \theta(x),\theta(y)\rangle \ =\ u\langle x,y\rangle  \quad (x,y\in E).
$$
\end{proposition}
\begin{proof}
Let
$u:=\langle \theta(x_0),\theta(x_0)\rangle  \in A_+$.
For any symmetry $w\in A$, as $x_0+x_0 w$ and $x_0 -x_0 w$ are orthogonal to each other,
so are $\theta(x_0 )+\theta(x_0)w$ and $\theta(x_0 )-\theta(x_0 )w$.
Consequently, $u + wu - uw - wuw = 0$ and $u + uw - wu - wuw = 0$ (by taking adjoint).
This tells us that $u = wuw$,
and so, $u\in Z(A)_+$ (as $A$ is generated by projections, and thus also by symmetries).
Pick any $z\in E$ with $\langle x_0,z\rangle  = 0$.
Then $z + x_0\langle z,z\rangle ^{1/2}$ is also orthogonal to $z - x_0\langle z,z\rangle ^{1/2}$.
It follows from the orthogonality preserving property that
$$
\langle \theta(z),\theta(z)\rangle
\ =\ \langle z,z\rangle ^{1/2} \langle \theta(x_0),\theta(x_0)\rangle  \langle z,z\rangle ^{1/2}
\ =\ u\langle z,z\rangle .
$$
For any $y\in E$, the element $z= y - x_0\langle x_0,y\rangle $ is orthogonal to $x_0$.
Hence,
$$\langle \theta(y),\theta(y)\rangle \
\ =\ \langle y,x_0\rangle \langle \theta(x_0),\theta(x_0)\rangle \langle x_0,y\rangle  +
\langle \theta(z),\theta(z)\rangle
\ =\ u \langle y,y\rangle .
$$
A polarization type argument implies that $\langle \theta(x),\theta(y)\rangle  = u\langle x,y\rangle $ ($x,y\in E$).
\end{proof}

\begin{theorem}\label{thm:unbdd-case-rro}
Let $A$ be a $C^*$-algebra of real rank zero.
Suppose that $E$ is full, and
$\theta: E\rightarrow F$ is an orthogonality preserving $A$-module map (not assumed to be bounded).
There is $u\in Z(M(A))_+$ such that
$$
\langle \theta(x),\theta(y)\rangle \ =\ u\langle x,y\rangle \quad (x,y\in E).
$$
In particular, $\theta$ is automatically bounded.
\end{theorem}
\begin{proof}
Set
$$P\ :=\ \{(x,p)\in E\times {\rm Proj}_1(A): \langle x, x \rangle = p \textrm{ and } xp = x\}.$$
Lemma \ref{lem:cp}(b) tells us that $P \neq \emptyset$.
Suppose that $(x,p)\in P$.
Then $E p$ is a full Hilbert $pAp$-module and the restriction of $\theta$ on $E p$ is an orthogonality preserving $pAp$-module map.
Since $p$ is the identity of the $C^*$-algebra $pAp$ and $\theta(Ep) \subseteq  Fp$, one can apply
Proposition \ref{prop:bddcase-with1} to obtain $b_p\in Z(pAp)_+$ that satisfies
$$
\langle \theta(x)p,\theta(y)p\rangle \ =\ b_p\langle xp,yp\rangle\qquad (x,y\in E).
$$
By Lemma \ref{lem:cp}(a), we have
\begin{equation}
\label{eqt:pap=0}
p\left(\langle \theta(x),\theta(y)\rangle  - \cv(b_p)\langle x,y\rangle \right)p\ =\ 0\qquad (x,y\in E).
\end{equation}
As the weak-*-closed linear span, $I$, of $\{\langle  \theta(x), \theta(y)\rangle  - \cv(b_p) \langle x, y\rangle: x,y\in E\}$ is an ideal of $A^{**}$, there is a central projection $q_I\in A^{**}$ with $I = q_IA^{**}$.
Since $pq_I = pq_Ip = 0$, we have $\cv(p)\leq 1- q_I$.
Consequently,
\begin{equation}
\label{eqt-cbp}
\cv(p) \langle \theta(x),\theta(y)\rangle \ = \ \cv(p)\cv(b_p)\langle x,y\rangle  \qquad (x,y\in E).
\end{equation}

Now, let
$\mathcal{D}\ :=\ \left\{ D \subseteq P: \cv(p)\cv(q) = 0 \textrm{ whenever }(x,p),(y,q)\in D\right\}$.
If we equip $\mathcal{D}$ with the usual inclusion, then Zorn's Lemma gives a maximal element $D_0 = \{(x_\gamma,p_\gamma)\}_{\gamma\in \Gamma}\in \mathcal{D}$.
Set $q_0 : =\bigvee_{\gamma\in \Gamma} \cv(p_\gamma)$ in ${\rm Proj}_1(A^{**})$, which is a central element.
Observes that $1 - q_0$ will not dominate a non-trivial open projection.
Indeed, if $0 \neq q\leq 1 - q_0$ is an open projection, then Lemma \ref{lem:cp}(b) produces an element $(y, r)\in P$ such that $r \leq q$.
Therefore, $D_0 \cup \{(y,r)\} \in \mathcal{D}$ which contradicts the maximality of $D_0$.
We now claim that the $*$-homomorphism $\Phi: A\rightarrow q_0A\subseteq A^{**}$ defined by $\Phi(a) = q_0a$ is injective.
Suppose on the contrary that there exists $a\in A_+$ with $\|a\| = 1$ and $\Phi(a) = 0$.
Take any $\epsilon \in (0,1)$, and put $q_\epsilon$ to be the non-zero open projection $e_a(\epsilon,1]$.
As $a - \epsilon q_\epsilon \geq 0$, we have $q_\epsilon q_0 = q_0q_\epsilon q_0\leq q_0aq_0/\epsilon =0$.
So, $q_\epsilon\leq 1-q_0$ which implies the contradiction that $q_\epsilon=0$.

As $x_{\gamma}\cv(p_\gamma) = x_{\gamma}p_\gamma\cv(p_\gamma) = x_\gamma$ ($\gamma\in \Gamma$), we see that $x_{\gamma}$ and $x_{\gamma'}$ are orthogonal if $\gamma\neq \gamma'$.
We now claim that $\cv(b_\gamma)$'s are uniformly bounded (where $b_\gamma\in Z(p_\gamma A p_\gamma)_+$
is the element associated with $(x_\gamma, p_\gamma)\in P$ that satisfies \eqref{eqt-cbp}).
Suppose on the contrary that there are $\cv(b_{\gamma_n})$ with $\|\cv(b_{\gamma_n})\|
= \|b_{\gamma_n}\|\geq n^3$ ($n\in \mathbb{N}$).
Note that the orthogonal sum $x:=\sum_n \frac{x_{\gamma_n}}{n}$ convergent in norm in $E$.
By the orthogonality preserving property of $\theta$, Lemma \ref{lem:cp}(a) as well as Equality \eqref{eqt:pap=0}, for any $m \in \mathbb{N}$,
\begin{align*}
\left< \theta(x),\theta(x)\right>
&=\  \left< \theta\left({x_{\gamma_m}/{m}}\right), \theta\left({x_{\gamma_m}/{m}}\right)\right>
+ \left< \theta\left(x - {x_{\gamma_m}/{m}}\right),\theta\left(x - {x_{\gamma_m}/{m}}\right)\right> \\
&\geq\ \left< \theta\left({x_{\gamma_m}/{m}}\right), \theta\left({x_{\gamma_m}/{m}}\right)\right>
\ = \  \frac{\cv(b_{\gamma_m})\left<  x_{\gamma_m},  x_{\gamma_m}\right> }{m^2}
\ = \  \frac{b_{\gamma_m}}{m^2}.
\end{align*}
As the norm of the last term goes to infinity as $n\to \infty$, we reach a contradiction.

Finally, let $d$ be the weak-*-limit in $A^{**}$ of finite
sums of the uniformly bounded mutually
orthogonal elements $\cv(b_\gamma)$ (see Lemma \ref{lem:cp}(a)).
By Relation \eqref{eqt-cbp} and the fact that $q_0$ is the weak-*-limit of finite sums of $\cv(p_\gamma)$'s, we have
$$d q_0 \langle x, y\rangle\ =\ q_0\langle \theta(x), \theta(y)\rangle\ \in\ q_0A \qquad (x,y\in E).$$
Since $E$ is full, we see that $d$ induces an element $m\in Z(M(q_0A))_+$ such that $m q_0 \langle x, y\rangle = q_0\langle \theta(x), \theta(y)\rangle$ ($x,y\in E$).
Since $\Phi:A\rightarrow q_0A$ extends to a $*$-isomorphism $\tilde \Phi: M(A) \rightarrow M(q_0A)$, there is $u\in Z(M(A))_+$ such that $\tilde \Phi(u) = m$.
This means that
$$\Phi(u \langle x, y\rangle - \langle \theta(x), \theta(y)\rangle) = 0 \qquad (x,y\in E)$$
which gives the required conclusion.
\end{proof}

\begin{remark}
Let $A$ be a general $C^*$-algebra. Suppose that there exist Hilbert
$A^{**}$-modules $\tilde E$ and $\tilde F$ containing  $E$ and $F$
respectively, such that the Hilbert $A^{**}$-module structures
extend the corresponding Hilbert $A$-module structures, and that
$\theta$ extends to an orthogonality preserving $A^{**}$-module map
$\tilde \theta: \tilde E \to \tilde F$. Then one can use Theorem
\ref{thm:unbdd-case-rro} to show that $\theta$ satisfies the
conclusion of Conjecture \ref{conj} (since $A^{**}$ has real rank
zero). In the situation when $\theta$ is a bounded orthogonality
preserving $A$-module map, we have tried $\tilde E = E^{**}$ and
$\tilde F = F^{**}$ but encountered some difficulties in showing
that $\theta^{**}$ is orthogonality preserving. 
It was claimed in \cite{FMP} that, when $\theta$ is bounded, such $\tilde E$,
$\tilde F$ and $\tilde \theta$ could  be found. 
However, instead of
manipulating the difficulties in the arguments in \cite{FMP}, we are
working on a proof of the general case, \emph{without} the
boundedness assumption on $\theta$, using completely
different ideas from those in this article, in \cite{FMP}, in
\cite{Turnsek-JMAA}, nor in \cite{LNW-orth-pres-comm}.
\end{remark}

\section{Orthogonality preserving $\mathbb C$-linear local maps}

In this section, we consider ($\mathbb{C}$-linear) local maps (see the Introduction) that preserve orthogonality.
Let us first give the following useful observation.

\begin{lemma}
\label{lem:local-A0-lin}
Let $A_0$ be the $*$-algebra generated by all the idempotents in $A$.
If $\Psi: E\rightarrow F$ is a local map, then $\Psi$ is an $A_0$-module map.
\end{lemma}
\begin{proof}
Let $p\in A$ be an idempotent and $x\in E$.
As $\Psi$ is local, one has $\Psi(x - xp)p = 0$.
If $\{u_i\}$ is an approximate unit for $A$, then $(1-p)u_i\in A$ will strictly converge to $(1-p)$.
Since $\Psi (xp) (1-p)u_i = 0$,
we have $\langle y, \Psi(xp)\rangle (1-p) = \lim\ \langle y, \Psi(xp)\rangle (1-p)u_i = 0$ ($y\in F$).
This implies that $\Psi(xp) - \Psi(xp)p = \Psi(xp)(1-p) = 0$.
Thus, $\Psi(x) p = \Psi(xp)$, and so $\Psi(xv)=\Psi(x)v$ for any $v\in A_0$.
\end{proof}

Note that if $A$ has real rank zero, then $A_0$ is dense in $A$.
We remark however that $A_0$ can be $\{0\}$ (e.g. if $A=C_0(0,1)$).

Recall that $A$ is a \emph{standard $C^*$-algebra} on a Hilbert space $H$ if $\mathcal{K}(H)\subseteq A\subseteq \CL(H)$.
In this case, $A_0$ contains a ``big enough''
ideal $\mathcal{F}(H)$ of $A$, in the sense that $\CK(H) = \overline{\mathcal{F}(H)}$ is an essential ideal. 
As a consequence, we can use Lemma \ref{lem:local-A0-lin} and Theorem
\ref{thm:unbdd-case-rro} to give a self-contained proof of the following slight extension of \cite[2.3]{Turnsek-JMAA} (note that the $A$-linearity is replaced by
the local property).

\begin{corollary}(c.f. \cite[2.3]{Turnsek-JMAA})
\label{cor:ortho-local-standard}
Suppose that $A$ is a standard $C^*$-algebra on a Hilbert space $H$.
If $\Psi: E\rightarrow F$ is local and orthogonality preserving, then
there is $\lambda\in \mathbb{R}_+$ such that
$$
\langle \Psi(x),\Psi(y)\rangle \ =\ \lambda\langle x,y\rangle  \qquad (x,y\in E).
$$
\end{corollary}
\begin{proof}
By Lemma \ref{lem:local-A0-lin}, $\Psi(xv)=\Psi(x)v$ for every finite rank operator $v$.
Consider $E_0:= E\cdot \CK(H)$ and $F_0:= F\cdot \CK(H)$.
Let $\{v_\gamma\}_{\gamma\in \Gamma}$ be an approximate unit in $\CK(H)$ consisting of finite rank positive operators.
For any $y\in E_0$, there exist $a\in \CK(H)$ and $x\in E_0$ with $y = xa$ (by the Cohen factorization theorem), and
$$\Psi(y)v_\gamma\ =\
\Psi(xa v_\gamma)\ =\ \Psi(x)a v_\gamma \quad (\gamma\in \Gamma)$$
which shows that $\|\Psi(y)v_\gamma - \Psi(x)a\| \to 0$ (along $\gamma$).
Define $\Phi: E_0\rightarrow F_0$ by setting $\Phi(y)$ to be the norm limit of $\Psi(y)v_\gamma$.
As $\Phi(y) = \Psi(x)a$ as well, we see that $\Phi(y)$ depends neither on the choice of $\{v_\gamma\}_{\gamma\in \Gamma}$, nor
on the decomposition $y = xa$.
If $b\in \CK(H)$, then
$$\Phi(yb)\ =\ \Phi(xab)\ =\ \Psi(x)ab\ =\ \Phi(y)b.$$
Moreover, if $x,y\in E_0$ with $\langle x, y\rangle =0$, then
$\langle \Psi(x), \Psi(y) \rangle = 0$ which implies that $\langle \Psi(x)v_\gamma, \Psi(y)v_{\gamma'} \rangle = 0$ ($\gamma, \gamma'\in \Gamma$), and so, $\langle \Phi(x), \Phi(y) \rangle = 0$.
On the other hand, since $\CK(H)$ is simple, we see that either $E_0$ is a full $\mathcal{K}(H)$-module or $E_0 = \{0\}$.
By Theorem \ref{thm:unbdd-case-rro}, there exists $\lambda\in Z(M(\CK(H)))_+ = \mathbb{R}_+$ such that for every $x,y\in E_0$, one has $\langle \Phi(x), \Phi(y)\rangle = \lambda \langle x, y\rangle$ (note that one can take any $\lambda$ if $E_0 = \{0\}$).
Thus, for any $x,y\in E$ and $\gamma, \gamma'\in \Gamma$,
$$
v_\gamma\langle \Psi(x),\Psi(y)\rangle v_{\gamma'}
\ =\ \langle \Phi(xv_\gamma), \Phi(yv_{\gamma'})\rangle
\ =\ \lambda\langle xv_\gamma,yv_{\gamma'}\rangle
\ =\ \lambda v_\gamma \langle x,y\rangle v_{\gamma'}.
$$
Consequently, if $b_{x,y}:= \langle \Psi(x),\Psi(y)\rangle  - \lambda \langle x,y\rangle  \in A$, then $v_\gamma b_{x,y} v_{\gamma'} = 0$ ($\gamma, \gamma'\in \Gamma$), which show that $b_{x,y} = 0$ (as $v_\gamma \to 1$ in the strong operator topology).
\end{proof}

\begin{corollary}
\label{cor:infiniteWalg}
Let $A$ be a  $W^*$-algebra with no finite type $\II$ direct summand.
If $E$ is full and $\Psi: E\rightarrow F$ is an orthogonality preserving local map, then there is $u\in Z(A)_+$ such that
$$
\langle \Psi(x),\Psi(y)\rangle \ =\ u\langle x,y\rangle\qquad (x,y\in E).
$$
\end{corollary}
\begin{proof}
By Theorem \ref{thm:unbdd-case-rro}, it suffices to show that $\Psi$ is an $A$-module map.
Recall that there are mutually orthogonal central projections
$q_{11}$, $q_{21}$ and $q_\infty$ in $A$ summing up to $1$ such that $q_{11} A$ is a finite $W^*$-algebra of type $I$, $q_{21} A$ is a finite $W^*$-algebra of type $\II$, and $q_\infty A$ is a properly infinite $W^*$-algebra (see, e.g., \cite[6.1.9]{Li}).
Since $A$ contains no finite type $\II$ direct summand, one has $q_{21} =0$.
Thus, $E = Eq_{11}\oplus Eq_\infty$.
The restriction $\Psi|_{Eq_\infty}$ is an $(q_\infty A)$-module map because of Lemma \ref{lem:local-A0-lin} and the fact that every element $a\in q_\infty A$ is a  sum of at most five idempotents
(see \cite[Theorem 4]{pearcy}).
Thus, it remains to verify the case when $A$ is a finite  $W^*$-algebra of type $I$.

In this case, for each $n\in \mathbb{N}$, there exist a hyperstonean space $\Omega_n$ (could be empty) and a projection $q_n\in Z(A)$ such that $\{q_n\}$ are orthogonal to one another, $\sum_n q_n$ weak-*-converges $1$ and $q_n A \cong C(\Omega_n)\otimes M_n$ (see e.g. \cite[6.7.7]{Li}).
Here we use the convention that $C(\Omega_n) = \{0\}$ if $\Omega_n = \emptyset$.
Let $n\in \mathbb{N}$ such that $\Omega_n \neq \emptyset$ and $e\in C(\Omega_n)$ be the identity.
Pick any rank one projection $p\in M_n$.
If $r:= e \otimes p\in {\rm Proj}_1(q_nA)$, then $rAr$ is isomorphic to $C(\Omega_n)$.
By Lemma \ref{lem:local-A0-lin}, the induced map $\Psi_r: Er\to Fr$ is an orthogonality preserving local map between Hilbert $rAr$-modules.
Using \cite[3.5]{LNW-orth-pres-comm}, we see that $\Psi_r$ is a $rAr$-module map.
In particular, for any $a\in C(\Omega_n)$ and $x\in E$, one has
$$
\Psi(x(a\otimes p))\ =\ \Psi_r(xr(a\otimes I)r)\ =\
\Psi_r(xr)r(a\otimes I)r\ =\ \Psi(x)(a\otimes p),
$$
where $I\in M_n$ is the identity.
Now, let $\{e_{kl}\}_{k,l=1}^n$ be the matrix unit of $M_n$.
As $\frac{e_{kl} + e_{kl}^* + e_{kk} + e_{ll}}{2}$ and $\frac{\mathrm{i}(e_{kl} - e_{kl}^*) + e_{kk} + e_{ll}}{2}$ are rank one projections, we see that $e_{kl}$ is a linear combinations of rank one projections.
Since any $a\in q_nA$ is of the form $a= \sum_{i,j = 1}^n a_{ij}\otimes e_{ij}$ ($a_{ij}\in C(\Omega_n)$), we see that $\Psi(xa) = \Psi(x)a$ ($x\in E; a\in q_nA$).
It follows that for any $x\in E$ and $a\in A$, 
$$\Psi(x a) \sum_{k=1}^n q_k\ =\ \Psi\left(x a \sum_{k=1}^n q_k\right)\ =\ \Psi(x) \left(a \sum_{k=1}^n q_k\right)\ =\ (\Psi(x)a) \sum_{k=1}^n q_k \qquad (n\in \mathbb{N}).$$
Consequently, for any $y\in F$, we have
$\langle y, \Psi(x a) - \Psi(x)a\rangle \sum_{k=1}^n q_k = 0$ which implies that $\langle y, \Psi(x a) - \Psi(x)a\rangle = 0$, and so $\Psi(x a) = \Psi(x)a$ as required.
\end{proof}

We end this article with a result that could be a first step towards a positive answer
for Conjecture \ref{conj} in the case when $A$ has real rank zero (with $\theta$ not being assumed to be an $A$-module map nor bounded).
This result is also interesting by itself, and gives us a rough idea what kind of difficulties will come across without the $A$-linearity.

\begin{proposition}\label{prop:local-rr0-with1}
Let $A$ be a unital $C^*$-algebra of real rank zero, and $A_0$ be the $*$-algebra generated by the idempotents in $A$.
Suppose that there is an element $x_0\in E$ such
that $\langle x_0,x_0\rangle  = 1$.
If $\Psi:E\to F$ is a local map preserving orthogonality, then one can find
$u\in Z(A)_+$ as well as an $A_0$-submodule $E_0\subseteq E$ containing $x_0$ with $E_0^\bot = \{0\}$ such that
$$
\langle \Psi(x),\Psi(y)\rangle \ =\ u\langle x,y\rangle  \quad (x,y\in E_0).
$$
\end{proposition}
\begin{proof}
Define $u:=\left<\Psi(x_0),\Psi(x_0)\right> \in A_+$.
Note that by Lemma \ref{lem:local-A0-lin}, $\Psi$ is an $A_0$-module map.
Thus, $\Psi(xw) = \Psi(x)w$ for any symmetry $w\in A$ (as $w\in A_0$).
Now, the argument of Proposition \ref{prop:bddcase-with1} tells us that $u\in Z(A)_+$.
Let $z\in E$ such that $\langle x_0,z\rangle  = 0$ and $\langle z,z\rangle  \in A_0$.
Then $z + x_0$ is also orthogonal to $z - x_0\langle z, z\rangle $.
It follows from the orthogonality preserving property and the $A_0$-linearity of $\Psi$ that
$$
\langle \Psi(z),\Psi(z)\rangle
\ =\ \langle \Psi(x_0),\Psi(x_0)\rangle  \langle z,z\rangle
\ =\ u\langle z,z\rangle .
$$
Let $\mathcal{D} := \{ D\subseteq E: x_0\in D; \langle x, x \rangle \in A_0 \textrm{ and } \langle x, y \rangle
=0 \textrm{ for any }x\neq y\in D\}$.
Take any maximal element $M$ in $\mathcal{D}$, and define $E_0$ to be the linear spans of $x\cdot a$ ($x\in M; a\in A_0$).
For any $y\in E_0$, we know that $\langle y, y\rangle, \langle x_0,y\rangle \in A_0$.
Thus, $z= y - x_0\langle x_0,y\rangle $ is orthogonal to $x_0$ and $\langle z,z\rangle  = \langle y,y\rangle  - \langle y,x_0\rangle \langle x_0,y\rangle \in A_0$.
Hence, the above implies that
$$\langle \Psi(y),\Psi(y)\rangle \
\ =\ \langle y,x_0\rangle \langle \Psi(x_0),\Psi(x_0)\rangle \langle x_0,y\rangle  +
\langle \Psi(z),\Psi(z)\rangle
\ =\ u \langle y,y\rangle .
$$
A polarization type argument tells us that $\langle \Psi(x),\Psi(y)\rangle  = u\langle x,y\rangle $ ($x,y\in E_0$).
Suppose on the contrary that there exists $z\in E$ with $\|z\| =1$ and $\langle z,x\rangle  =0$ for any $x\in E_0$.
Let $a:=\langle z,z\rangle $ and $q_n := e_a(\frac{1}{2^n},1]$ ($n\in \mathbb{N}$).
There exist $d,b\in C^*(a)_+$ such that $q_5\leq ab\leq 1$, $d \leq a$, $dq_4 = aq_4$ and $d q_{5} = d$.
As $b^{1/2}d^{1/2} \leq 1$, $b^{1/2}d^{1/2}q_4 = b^{1/2}a^{1/2}q_4 = q_4$ and $b^{1/2}d^{1/2}q_5 = b^{1/2}d^{1/2}$, we see that
$$
\|z - zb^{1/2}d^{1/2}\|^2
\ = \ \|a - 2ab^{1/2}d^{1/2} + abd\|
\ \leq \ 2 \|a(1 - b^{1/2}d^{1/2})\|
\ < \ 1/8.
$$
Since $d^{1/2}\in \her(q_{5})$, there exists $c \in A_0\cap \her(q_{5})_+$ such that $\|d^{1/2} - c \| < 1/8$ (because $\her(q_{5})$ also has real rank zero; see e.g. \cite{Brown91}).
If $z' := zb^{1/2}c$, then $\langle z', z' \rangle = cq_{5}abq_{5}c = c^2\in A_0$ (as $abq_{5} = q_{5}$) and
$\|z - z'\| \leq \|z- zb^{1/2}d^{1/2}\| +\| zb^{1/2}d^{1/2} - zb^{1/2}c\| \leq 1/2$, which implies that $z' \neq 0$.
Moreover, $\langle x, z'\rangle = \langle x, z \rangle b^{1/2}c = 0$ for any $x\in M$, we see that $M\cup \{z'\} \in \mathcal{D}$, which is a contradiction.
\end{proof}


\begin{thebibliography}{99}

\bibitem{Al08}
J. Alaminos, M. Bre\v{s}ar, M. \v{C}erne, J. Extremera, A. R. Villena,
Zero product preserving maps on $C^1[0,1]$,
J. Math. Anal. Appl. \textbf{347} (2008), 472--481.

\bibitem{Ar04}
J. Araujo, Linear biseparating maps between spaces of
vector-valued differentiable functions and automatic continuity,
Adv. Math. \textbf{187} (2004), no. 2, 488--520.

\bibitem{Blanco06}
A. Blanco and A. Turn\v{s}ek, On maps that preserve orthogonality in normed spaces,
Proc. Royal Soc. Edingburgh \textbf{136A} (2006), 709--716.

\bibitem{Bro88} L. G. Brown, Semicontinuity and multipliers of $C^*$-algebras, Can.\  J. Math., \textbf{vol.\ {\rm XL}}, no.\ 4 (1988), 865--988.

\bibitem{Brown91}
L. G. Brown and G. K. Pedersen, $C^*$-algebras of real rank
zero, J. Funct.\ Anal.\ \textbf{99} (1991), 131--149.

\bibitem{chmielinski05}
J. Chmieli\'{n}ski, Linear mappings approximately preserving orthogonality,
J. Math. Anal. Appl. \textbf{304} (2005), 158--169.

\bibitem{FMP}
M. Frank, A.S. Mishchenko and A.A. Pavlov, Orthogonality-preserving, $C^*$-conformal and conformal module mappings on Hilbert $C^*$-module, preprint (arXiv:0907.2983).

\bibitem{Turnsek-JMAA}
D. Ili\v{s}evi\'{c} and A. Turn\v{s}ek, Approximately orthogonality preserving mappings on $C^*$-modules,
J. Math. Anal. Appl. \text{341} (2008), 298--308.

\bibitem{KR}
R. V. Kadison and J. R. Ringrose, \emph{Fundamentals of the theory of operator algebras,
Vol.\ I}, Graduate Studies in Math, 15. Amer. Math. Soc., Providence, RI, 1997.

\bibitem{KN} R. Kantrowitz and M.M. Neumann, Disjointness preserving and local operators on algebras of differentiable functions, Glasg. Math. J. \textbf{43} (2001), 295--309.

\bibitem{LNW-auto-cont} C.W. Leung, C.K. Ng and N.C. Wong, Automatic continuity and $C_0(\Omega)$-linearity
of linear maps between $C_0(\Omega)$-modules, preprint.

\bibitem{LNW-orth-pres-comm}
C.W. Leung, C.K. Ng and N.C. Wong, Linear orthogonality
preservers of Hilbert bundles, preprint.

\bibitem{Li}
Bingren Li, \emph{Introduction to operator algebras}, World Scientific, Singapore, 1992.

\bibitem{Na85}
R. Narasimhan, \emph{Analysis on real and complex manifolds},
Advanced Studies in Pure Mathematics \textbf{1}, North-Holland
Publishing Co., Amsterdam (1968).

\bibitem{pearcy}
C. Pearcy and D. Topping, Sums of small numbers of
idempotents, Michigan Math. J. \textbf{14} (1967), 453--465.


\bibitem{Ped} G.K. Pedersen, \emph{$C^*$-algebras and their automorphism groups}, Academic Press, New York, (1979).

\bibitem{Pe60}
J. Peetre, R\'{e}ctification \`{a} l'article ``Une
caract\'{e}risation abstraite des op\'{e}rateurs
diff\'{e}rentiels'', Math. Scand. \textbf{8} (1960), 116--120.

\end{thebibliography}
\end{document}